\documentclass[12pt,reqno]{article}

\usepackage[usenames]{color}
\usepackage{amssymb}
\usepackage{graphicx}
\usepackage{amscd}

\usepackage{amsthm}
\newtheorem{theorem}{Theorem}

\newtheorem{proposition}[theorem]{Proposition}
\newtheorem{corollary}[theorem]{Corollary}

\theoremstyle{definition}

\newtheorem{example}[theorem]{Example}

\usepackage[colorlinks=true,
linkcolor=webgreen, filecolor=webbrown,
citecolor=webgreen]{hyperref}

\definecolor{webgreen}{rgb}{0,.5,0}
\definecolor{webbrown}{rgb}{.6,0,0}

\usepackage{color}

\usepackage{float}

\usepackage{graphics,amsmath,amssymb}
\usepackage{amsfonts}
\usepackage{latexsym}
\usepackage{epsf}

\newcommand{\seqnum}[1]{\href{http://www.research.att.com/cgi-bin/access.cgi/as/~njas/sequences/eisA.cgi?Anum=#1}{\underline{#1}}}

\begin{document}

\begin{center}
\vskip 1cm{\LARGE\bf Simple closed form Hankel transforms based on
the central coefficients of certain Pascal-like triangles} \vskip
1cm \large
Paul Barry\\
School of Science\\
Waterford Institute of Technology\\
Ireland\\
\href{mailto:pbarry@wit.ie}{\tt pbarry@wit.ie} \\
\end{center}
\vskip .2 in

\begin{abstract} We study the Hankel transforms of sequences related to the central
coefficients of a family of Pascal-like triangles. The mechanism of
Riordan arrays is used to elucidate the structure of these
transforms.
\end{abstract}
\section{Introduction}
This note concerns the characterization of the Hankel transfoms of
the central coefficients $T(2n,n,r)$ of a family of Pascal-like
triangles that are parameterised by an integer $r$. Specifically, we
define a family of number triangles with general term $T(n,k,r)$ by
$$T(n,k,r)=\sum_{k=0}^{n-k}\binom{k}{j}\binom{n-k}{j}r^j.$$ For
instance, $r=1$ gives Pascal's triangle \seqnum{A007318}, while
$r=2$ gives the triangle of Delannoy numbers \cite{PasTri},
\seqnum{A008288}.
\begin{proposition} The Hankel transform of the sequence
$a(n,r)=T(2n,n,r)$ is given by
$$2^n r^{\binom{n+1}{2}}.$$ \end{proposition}
\begin{proof} We proceed as in \cite{WW} and \cite{PPWW} by means of
the $LDL^T$ decomposition of the Hankel matrix $H(r)$ of
$T(2n,n,r)$. We take the example of $r=2$. In this case,
$$H(2)=\left(\begin{array}{ccccc} 1 & 3 & 13
& 63 & \ldots \\3 & 13 & 63 & 321 & \ldots \\ 13 & 63 & 321 & 1683 &
\ldots \\ 63 & 321 & 1683 & 8989 & \ldots \\  \vdots & \vdots &
\vdots & \vdots & \ddots\end{array}\right)$$ Then
\begin{eqnarray*}H(2)&=&L(2)D(2)L(2)^T\\
&=&\left(\begin{array}{ccccc} 1 & 0 & 0 & 0 & \ldots \\3 & 1 & 0 & 0
& \ldots \\ 13 & 6 & 1 & 0 & \ldots \\ 63 & 33 & 9 & 1 & \ldots
\\  \vdots & \vdots & \vdots & \vdots &
\ddots\end{array}\right)\left(\begin{array}{ccccc} 1 & 0 & 0 & 0 &
\ldots \\0 & 4 & 0 & 0 & \ldots \\ 0 & 0 & 8 & 0 & \ldots \\ 0 & 0 &
0 & 16 & \ldots
\\  \vdots & \vdots & \vdots & \vdots &
\ddots\end{array}\right)\left(\begin{array}{ccccc} 1 & 3 & 13 & 63 &
\ldots \\0 & 1 & 6 & 33 & \ldots \\ 0 & 0 & 1 & 9 & \ldots \\ 0 & 0
& 0 & 1 & \ldots
\\  \vdots & \vdots & \vdots & \vdots &
\ddots\end{array}\right)\end{eqnarray*} Hence the Hankel transform
of $T(2n,n,2)$ is equal to the sequence with general term
$$\prod_{k=0}^n(2.2^k-0^k)=2^n 2^{\binom{n+1}{2}}.$$ $L(2)$ is in
fact the Riordan array
$$(\frac{1}{\sqrt{1-6x+x^2}},\frac{1-3x-\sqrt{1-6x+x^2}}{4x})$$ or
$$(\frac{1-2x^2}{1+3x+2x^2},\frac{x}{1+3x+2x^2})^{-1}.$$
In general, we can show that $H(r)=L(r)D(r)L(r)^T$ where $L(r)$ is
the Riordan array
$$(\frac{1}{\sqrt{1-2(r+1)x+(r-1)^2x^2}},\frac{1-(r+1)x-\sqrt{1-2(r+1)x+(r-1)^2x^2}}{2rx})$$
and $D(r)$ is the diagonal matrix with $n$-th term $2.r^n-0^n$.
Hence the Hankel transform of $T(2n,n,r)$ is given by
$$\prod_{k=0}^n(2.r^k-0^k)=2^n r^{\binom{n+1}{2}}.$$
\end{proof}
We note that the Riordan array $L(r)$
$$(\frac{1}{\sqrt{1-2(r+1)x+(r-1)^2x^2}},\frac{1-(r+1)x-\sqrt{1-2(r+1)x+(r-1)^2x^2}}{2rx})$$
is the inverse of the Riordan array
$$(\frac{1-rx^2}{1+(r+1)x+rx^2},\frac{x}{1+(r+1)x+rx^2}).$$
Its general term is given by
$$\sum_{j=0}^n\binom{n}{j}\binom{n}{j-k}r^{j-k}=\sum_{j=0}^n
\binom{n}{j}\binom{j}{n-k-j}r^{n-k-j}(r+1)^{2j-(n-k)}.$$ Its $k$-th
column has exponential generating function given by
$$e^{(r+1)x}I_k(2\sqrt{r}x)/{\sqrt{r}}^k.$$
\begin{corollary} The sequences with e.g.f. $I_0(2\sqrt{r}x)$ have
Hankel transforms given by $2^nr^{\binom{n+1}{2}}$. \end{corollary}
\begin{proof} By \cite{PasTri} or otherwise, we know that the
sequences $T(2n,n,r)$ have e.g.f.
$$e^{(r+1)x}I_0(2\sqrt{r}x).$$
By the above proposition and the binomial invariance property
 of the Hankel transform \cite{Layman}, $\mathbf{B}^{-r-1}T(2n,n,r)$
 has the desired Hankel transform. But $\mathbf{B}^{-r-1}T(2n,n,r)$
 has e.g.f. given by
 $$e^{-(r+1)x}e^{(r+1)x}I_0(2\sqrt{r}x)=I_0(2\sqrt{r}x).$$ \end{proof}
\section{Hankel transform of generalized Catalan numbers}
Following \cite{PasTri}, we denote by $c(n;r)$ the sequence of
numbers
$$ c(n;r)=T(2n,n,r)-T(2n,n+1,r).$$
For instance, $c(n;1)=c(n)$, the sequence of Catalan numbers
\seqnum{A000108}. We have
\begin{proposition} The Hankel transform of $c(n;r)$ is
$r^{\binom{n+1}{2}}$. \end{proposition} \begin{proof} Again, we use
the $LDL^T$ decomposition of the associated Hankel matrices. For
instance, when $r=3$, we obtain
$$H(3)=\left(\begin{array}{ccccc} 1 & 3 & 12
& 57 & \ldots \\3 & 12 & 57 & 300 & \ldots \\ 12 & 57 & 300 & 1686 &
\ldots \\ 57 & 300 & 1686 & 9912 & \ldots \\  \vdots & \vdots &
\vdots & \vdots & \ddots\end{array}\right)$$ Then
\begin{eqnarray*}H(3)&=&L(3)D(3)L(3)^T\\
&=&\left(\begin{array}{ccccc} 1 & 0 & 0 & 0 & \ldots \\3 & 1 & 0 & 0
& \ldots \\ 12 & 7 & 1 & 0 & \ldots \\ 57 & 43 & 11 & 1 & \ldots
\\  \vdots & \vdots & \vdots & \vdots &
\ddots\end{array}\right)\left(\begin{array}{ccccc} 1 & 0 & 0 & 0 &
\ldots \\0 & 3 & 0 & 0 & \ldots \\ 0 & 0 & 9 & 0 & \ldots \\ 0 & 0 &
0 & 27 & \ldots
\\  \vdots & \vdots & \vdots & \vdots &
\ddots\end{array}\right)\left(\begin{array}{ccccc} 1 & 3 & 12 & 57 &
\ldots \\0 & 1 & 7 & 43 & \ldots \\ 0 & 0 & 1 & 11 & \ldots \\ 0 & 0
& 0 & 1 & \ldots
\\  \vdots & \vdots & \vdots & \vdots &
\ddots\end{array}\right)\end{eqnarray*} Hence the Hankel transform
of $c(n;3)$ is $$\prod_{k=0}^n3^k=3^{\binom{n+1}{2}}.$$ In this
case, $L(3)$ is the Riordan array
$$(\frac{1}{1+3x},\frac{x}{1+4x+3x^2})^{-1}.$$
In general, we can show that $H(r)=L(r)D(r)L(r)^T$ where
$$L(r)=(\frac{1}{1+rx},\frac{x}{1+(r+1)x+rx^2})^{-1}$$ and $D(r)$ has $n$-th term
$r^n$. Hence the Hankel transform of $c(n;r)$ is given by
$$\prod_{k=0}^nr^k=r^{\binom{n+1}{2}}.$$
\end{proof}
We finish this section with some notes concerning production
matrices as found, for instance, in \cite{DFR}. It is well known
that the production matrix $P(1)$ for the Catalan numbers
$C(n)=c(n,1)$ is given by
$$P(1)=\left(\begin{array}{ccccc} 0 & 1 & 0 & 0 & \ldots \\0 & 1 & 1 & 0 &
\ldots \\ 0 & 1 & 1 & 1 & \ldots
\\  \vdots & \vdots & \vdots & \vdots &
\ddots\end{array}\right)$$ Following \cite{DFR}, we can associate a
Riordan array $A_P(1)$ to $P(1)$ as follows. The second column of
$P$ has generating function $\frac{1}{1-x}$. Solving the equation
$$u=\frac{1}{1-xu}$$ we obtain $u(x)=\frac{1-\sqrt{1-4x}}{2x}=c(x)$.
Since the first column is all $0$'s, this means that $A_P(1)$ is the
Riordan array $(1,xc(x))$. This is the inverse of $(1,x(1-x))$. We
have $$A_P(1)=\left(\begin{array}{ccccc} 1 & 0 & 0 & 0 & \ldots \\0
& 1 & 0 & 0 & \ldots \\ 0 & 1 & 1 & 0 & \ldots
\\ 0 & 2 & 2 & 1 & \ldots \\  \vdots & \vdots & \vdots &
\vdots & \ddots\end{array}\right)$$ Multiplying on the right by $B$,
the binomial matrix, we obtain
$$A_P(1)B=\left(\begin{array}{ccccc} 1 & 0 & 0 & 0 & \ldots \\1
& 1 & 0 & 0 & \ldots \\ 2 & 3 & 1 & 0 & \ldots
\\ 5 & 9 & 5 & 1 & \ldots \\  \vdots & \vdots & \vdots &
\vdots & \ddots\end{array}\right)=L(1)$$ which is the Riordan array
$$(\frac{1}{1-x},xc(x)^2)=(\frac{1}{1+x},\frac{1}{1+2x+x^2})^{-1}.$$

 Similarly the production matrix
for the $c(n;2)$, or the large Schroeder numbers, is given by
$$P(2)=\left(\begin{array}{ccccc} 0 & 2 & 0 & 0 & \ldots \\0 & 1 & 2 & 0 &
\ldots \\ 0 & 1 & 1 & 2 & \ldots
\\  \vdots & \vdots & \vdots & \vdots &
\ddots\end{array}\right)$$ Here, the generating function for the
second column is $\frac{2-x}{1-x}$. Now solving
$$u=\frac{2-xu}{1-xu}$$ which gives
$u=\frac{1+x-\sqrt{1-6x+x^2}}{2x}$. Hence in this case, $A_P(2)$ is
the Riordan array $(1,\frac{1+x-\sqrt{1-6x+x^2}}{2})$. That is,
$$A_P(2)=\left(\begin{array}{ccccc} 1 & 0 & 0
& 0 & \ldots \\0 & 2 & 0 & 0 & \ldots \\ 0 & 2 & 4 & 0 & \ldots
\\ 0 & 6 & 8 & 8 & \ldots \\  \vdots & \vdots & \vdots &
\vdots & \ddots\end{array}\right)=(1,\frac{x(1-x)}{2-x})^{-1}.$$ The
row sums of this matrix are $1,2,6,22,90,\ldots$ as expected.
Multiplying $A_P(2)$ on the right by the binomial matrix $B$, we
obtain
$$A_P(2)B=\left(\begin{array}{ccccc} 1 & 0 & 0
& 0 & \ldots \\2 & 2 & 0 & 0 & \ldots \\ 6 & 10 & 4 & 0 & \ldots
\\ 22 & 46 & 32 & 8 & \ldots \\  \vdots & \vdots & \vdots &
\vdots & \ddots\end{array}\right)$$ which is the array
$$(\frac{1-x-\sqrt{1-6x+x^2}}{2x},\frac{1-3x-\sqrt{1-6x+x^2}}{2x}).$$ Finally
$$A_PB\left(\begin{array}{ccccc} 1 & 0 & 0
& 0 & \ldots \\0 & \frac{1}{2} & 0 & 0 & \ldots \\ 0 & 0 &
\frac{1}{4} & 0 & \ldots
\\ 0 & 0 & 0 & \frac{1}{8} & \ldots \\  \vdots & \vdots & \vdots &
\vdots &
\ddots\end{array}\right)=A_PB(1,\frac{x}{2})=\left(\begin{array}{ccccc}
1 & 0 & 0 & 0 & \ldots \\2 & 1 & 0 & 0 & \ldots \\ 6 & 5 & 1 & 0 &
\ldots
\\ 22 & 23 & 8 & 1 & \ldots \\  \vdots & \vdots & \vdots &
\vdots & \ddots\end{array}\right)=L(2)$$ which is
$$(\frac{1-x-\sqrt{1-6x+x^2}}{2x},\frac{1-3x-\sqrt{1-6x+x^2}}{4x})$$
or $$L(2)=(\frac{1}{1+2x},\frac{x}{1+3x+2x^2})^{-1}.$$ We can
generalize these results to give the following proposition.
\begin{proposition} The production matrix for the generalized
Catalan sequence $c(n;r)$ is given by
$$P(r)=\left(\begin{array}{ccccc} 0 & r & 0 & 0 & \ldots \\0 & 1 & r & 0 &
\ldots \\ 0 & 1 & 1 & r & \ldots
\\  \vdots & \vdots & \vdots & \vdots &
\ddots\end{array}\right)$$ The associated matrix $A_P(r)$ is given
by
$$A_P(r)=(1,\frac{x(1-x)}{r-(r-1)x})^{-1}=(1,\frac{1+(r-1)x-\sqrt{1-2(r+1)x+(r-1)^2x^2}}{2}).$$
The matrix $L(r)$ in the decomposition $L(r)D(r)L(r)^T$ of the
Hankel matrix $H(r)$ for $c(n;r)$, which is equal to
$A_P(r)B(1,x/r)$,  is given by
$$L(r)=(\frac{1-(r-1)x-\sqrt{1-2(r+1)x+(r-1)^2x^2}}{2x},\frac{1-(r+1)x-\sqrt{1-2(r+1)x+(r-1)^2x^2}}{2rx}).$$
We have $$L(r)=(\frac{1}{1+rx},\frac{x}{1+(r+1)x+rx^2})^{-1}.$$
\end{proposition} We note that the elements of $L(r)^{-1}$ are in fact the
coefficients of the orthogonal polynomials associated to $H(r)$.
\begin{proposition}The elements of the rows of the Riordan array
$(\frac{1}{1+rx},\frac{x}{1+(r+1)x+rx^2})$ are the coefficients of
the orthogonal polynomials associated to the Hankel matrix
determined by the generalized Catalan numbers
$c(n;r)$.\end{proposition}

\section{Hankel transform of the sum of consecutive generalized Catalan numbers}
We now look at the Hankel transform of the sum of two consecutive generalized Catalan numbers.
That is, we study the Hankel transform of $c(n;r)+c(n+1;r)$. For the case $r=1$ (the ordinary Catalan numbers) this
was dealt with in \cite{Hankel1}, while the general case was studied in \cite{Hankel2}. We use the methods
developed above to gain greater insight. We start with the case $r=1$. For this, the Hankel matrix for
$c(n)+c(n+1)$ is given by
$$H=\left(\begin{array}{ccccc} 2 & 3 & 7
& 19 & \ldots \\3 & 7 & 19 & 56 & \ldots \\ 7 & 19 & 56 & 174 &
\ldots \\ 19 & 56 & 174 & 561 & \ldots \\  \vdots & \vdots &
\vdots & \vdots & \ddots\end{array}\right)$$ Proceeding to the $LDL^T$ decomposition, we get
\begin{eqnarray*}H&=&LDL^T\\
&=&\left(\begin{array}{ccccc} 1 & 0 & 0 & 0 & \ldots \\\frac{3}{2} & 1 & 0 & 0
& \ldots \\ \frac{7}{2} & \frac{17}{5} & 1 & 0 & \ldots \\ \frac{19}{2} & 11 & \frac{70}{13} & 1 & \ldots
\\  \vdots & \vdots & \vdots & \vdots &
\ddots\end{array}\right)\left(\begin{array}{ccccc} 2 & 0 & 0 & 0 &
\ldots \\0 & \frac{5}{2} & 0 & 0 & \ldots \\ 0 & 0 & \frac{13}{5} & 0 & \ldots \\ 0 & 0 &
0 & \frac{34}{13} & \ldots
\\  \vdots & \vdots & \vdots & \vdots &
\ddots\end{array}\right)\left(\begin{array}{ccccc} 1 & \frac{3}{2} & \frac{7}{2} & \frac{19}{2} &
\ldots \\0 & 1 & \frac{17}{5} & 11 & \ldots \\ 0 & 0 & 1 & \frac{70}{13} & \ldots \\ 0 & 0
& 0 & 1 & \ldots
\\  \vdots & \vdots & \vdots & \vdots &
\ddots\end{array}\right)\end{eqnarray*}
This indicates that the Hankel transform of $c(n)+c(n+1)$ is given by
$$\prod_{k=0}^n \frac{F(2k+3)}{F(2k+1)}=F(2n+3).$$ This is in agreement with \cite{Hankel1}. We note that in this case,
$L^{-1}$ takes the form
$$L^{-1}=\left(\begin{array}{ccccc} 1 & 0 & 0 & 0 & \ldots \\-\frac{3}{2} & 1 & 0 & 0
& \ldots \\ \frac{8}{5} & -\frac{17}{5} & 1 & 0 & \ldots \\ -\frac{21}{13} & \frac{95}{13} & -\frac{70}{13} & 1 & \ldots
\\  \vdots & \vdots & \vdots & \vdots &
\ddots\end{array}\right)=\left(\begin{array}{ccccc} 1 & 0 & 0 & 0 & \ldots \\0 & \frac{1}{2} & 0 & 0
& \ldots \\ 0 & 0 & \frac{1}{5} & 0 & \ldots \\ 0 & 0 & 0 & \frac{1}{13} & \ldots
\\  \vdots & \vdots & \vdots & \vdots &
\ddots\end{array}\right)\left(\begin{array}{ccccc} 1 & 0 & 0 & 0 & \ldots \\-3 & 2 & 0 & 0
& \ldots \\ 8 & -17 & 5 & 0 & \ldots \\ -21 & 95 & -70 & 13 & \ldots
\\  \vdots & \vdots & \vdots & \vdots &
\ddots\end{array}\right)$$ where we see the sequences $F(2n+1)$ and $(-1)^nF(2n+2)$ in evidence.

Now looking at the case $r=2$, we get
$$H=\left(\begin{array}{ccccc} 3 & 8 & 28
& 112 & \ldots \\8 & 28 & 112 & 484 & \ldots \\ 28 & 112 & 484 & 2200 &
\ldots \\ 112 & 484 & 2200 & 10364 & \ldots \\  \vdots & \vdots &
\vdots & \vdots & \ddots\end{array}\right)$$ Proceeding to the $LDL^T$ decomposition, we obtain
\begin{eqnarray*}H&=&LDL^T\\
&=&\left(\begin{array}{ccccc} 1 & 0 & 0 & 0 & \ldots \\\frac{8}{3} & 1 & 0 & 0
& \ldots \\ \frac{28}{3} & \frac{28}{5} & 1 & 0 & \ldots \\ \frac{112}{3} & \frac{139}{5} & \frac{146}{17} & 1 & \ldots
\\  \vdots & \vdots & \vdots & \vdots &
\ddots\end{array}\right)\left(\begin{array}{ccccc} 3 & 0 & 0 & 0 &
\ldots \\0 & \frac{20}{3} & 0 & 0 & \ldots \\ 0 & 0 & \frac{272}{20} & 0 & \ldots \\ 0 & 0 &
0 & \frac{7424}{272} & \ldots
\\  \vdots & \vdots & \vdots & \vdots &
\ddots\end{array}\right)\left(\begin{array}{ccccc} 1 & \frac{8}{2} & \frac{28}{3} & \frac{112}{3} &
\ldots \\0 & 1 & \frac{28}{5} & \frac{139}{5} & \ldots \\ 0 & 0 & 1 & \frac{146}{17} & \ldots \\ 0 & 0
& 0 & 1 & \ldots
\\  \vdots & \vdots & \vdots & \vdots &
\ddots\end{array}\right)\end{eqnarray*} Thus the Hankel transform of
$c(n;2)+c(n+1;2)$ is $3,20,272,7424\ldots$. This is in agreement
with \cite{Hankel2}. We note that different factorizations of
$L^{-1}$ can lead to different formulas for $h_n(2)$, the Hankel
transform of $c(n;2)+c(n+1;2)$. For instance, we can show that
$$L^{-1}=\left(\begin{array}{ccccc} 1 & 0 & 0 & 0 & \ldots \\-\frac{8}{3} & 1 & 0 & 0
& \ldots \\ \frac{28}{5} & -\frac{28}{5} & 1 & 0 & \ldots \\ -\frac{192}{17} & \frac{345}{17} & -\frac{146}{17} & 1 & \ldots
\\  \vdots & \vdots & \vdots & \vdots &
\ddots\end{array}\right)=\left(\begin{array}{ccccc} 1 & 0 & 0 & 0 & \ldots \\0 & \frac{1}{3} & 0 & 0
& \ldots \\ 0 & 0 & \frac{1}{5} & 0 & \ldots \\ 0 & 0 & 0 & \frac{1}{17} & \ldots
\\  \vdots & \vdots & \vdots & \vdots &
\ddots\end{array}\right)\left(\begin{array}{ccccc} 1 & 0 & 0 & 0 & \ldots \\-8 & 3 & 0 & 0
& \ldots \\ 28 & -28 & 5 & 0 & \ldots \\ -192 & 345 & -146 & 17 & \ldots
\\  \vdots & \vdots & \vdots & \vdots &
\ddots\end{array}\right)$$ We note that the diagonal elements of the last matrix correspond to the sequence $a(n)$ of terms $1,3,5,17,29,99,\ldots$ with
generating function $$\frac{1+3x-x^2-x^3}{1-6x^2+x^4}.$$
This is \seqnum{A079496}. It is the interleaving of bisections of the Pell numbers \seqnum{A000129} and their associated numbers
\seqnum{A001333}. We have
\begin{eqnarray*}a(n)&=&\sum_{k=0}^{\lfloor \frac{n+1}{2} \rfloor}\binom{n+1}{2k}2^{n+1-k-{\lfloor \frac{n+2}{2} \rfloor}}\\
&=&-(\sqrt{2}-1)^n((\frac{\sqrt{2}}{8}-\frac{1}{4})(-1)^n-\frac{\sqrt{2}}{8}-\frac{1}{4})-(\sqrt{2}+1)^n((\frac{\sqrt{2}}{8}-\frac{1}{4})(-1)^n-\frac{\sqrt{2}}{8}-1/4)\end{eqnarray*}
Multiplying $a(n)$ by $4^{\lfloor \frac{(n+1)^2}{4}\rfloor}$, we obtain $1,3,20,272,7424,\ldots$.
Hence
\begin{eqnarray*}1,3,20,272,\ldots&=&4^{\lfloor \frac{(n+1)^2}{4}\rfloor}\sum_{k=0}^{\lfloor \frac{n+1}{2} \rfloor}\binom{n+1}{2k}2^{n+1-k-{\lfloor \frac{n+2}{2} \rfloor}}\\
&=&4^{\lfloor \frac{(n+1)^2}{4}\rfloor}2^{n+1-{\lfloor \frac{n+2}{2} \rfloor}}\sum_{k=0}^{\lfloor \frac{n+1}{2} \rfloor}\binom{n+1}{2k}2^{-k}\\
&=&2^{\binom{n+1}{2}}\sum_{k=0}^{\lfloor \frac{n+1}{2}
\rfloor}\binom{n+1}{2k}2^{-k}\end{eqnarray*} That is, the Hankel
transform $h_n(2)$ of $c(n;2)+c(n+1;2)$ is given by $$
h_n(2)=2^{\binom{n+2}{2}}\sum_{k=0}^{\lfloor \frac{n+2}{2}
\rfloor}\binom{n+2}{2k}2^{-k}.$$

For our purposes, the following factorization of $L^{-1}$ is more
convenient.
$$L^{-1}=\left(\begin{array}{ccccc} 1 & 0 & 0 & 0 & \ldots \\-\frac{8}{3} & 1 & 0 & 0
& \ldots \\ \frac{56}{10} & -\frac{56}{10} & 1 & 0 & \ldots \\
-\frac{384}{34} & \frac{690}{34} & -\frac{292}{34} & 1 & \ldots
\\  \vdots & \vdots & \vdots & \vdots &
\ddots\end{array}\right)=\left(\begin{array}{ccccc} 1 & 0 & 0 & 0 &
\ldots \\0 & \frac{1}{3} & 0 & 0 & \ldots \\ 0 & 0 & \frac{1}{10} &
0 & \ldots \\ 0 & 0 & 0 & \frac{1}{34} & \ldots
\\  \vdots & \vdots & \vdots & \vdots &
\ddots\end{array}\right)\left(\begin{array}{ccccc} 1 & 0 & 0 & 0 &
\ldots \\-8 & 3 & 0 & 0 & \ldots \\ 56 & -56 & 10 & 0 & \ldots \\
-384 & 690 & -292 & 34 & \ldots
\\  \vdots & \vdots & \vdots & \vdots &
\ddots\end{array}\right)$$

We now note that the sequence $\frac{h_n(2)}{2^{\binom{n+1}{2}}}$ is
the sequence $b_2(n+1)$, where $b_2(n)$ is the sequence $1, 3, 10,
34, 116, \ldots$ with generating function $\frac{1-x}{1-4x+2x^2}$
and general term $$ b_2(n)=\sum_{k=0}^{\lfloor \frac{n}{2}
\rfloor}\binom{n-k}{k}(-2)^k4^{n-2k}-\sum_{k=0}^{\lfloor
\frac{n-1}{2} \rfloor}\binom{n-k-1}{k}(-2)^k4^{n-2k-1}.$$ Hence
$$h_n(2)=2^{\binom{n+1}{2}}b_2(n+1).$$
Noting that $b_2(n)$ is the binomial transform of the Pell
\seqnum{A000129}$(n+1)$ numbers whose generating function is
$\frac{1}{1-2x-x^2}$, we have the following alternative expressions
for $b_2(n)$:
\begin{eqnarray*}b_2(n)&=&\sum_{k=0}^n\binom{n}{k}\sum_{j=0}^k
\binom{j}{k-j}2^{2j-k}\\
&=&\sum_{k=0}^n\binom{n}{k}\sum_{j=0}^{\lfloor \frac{k}{2} \rfloor}
\binom{k-j}{j}2^{k-2j}.\end{eqnarray*}
 For $r=3$, we have
$$H=\left(\begin{array}{ccccc} 4 & 15 & 69
& 357 & \ldots \\15 & 69 & 357 & 1986 & \ldots \\ 69 & 357 & 1986 &
11598 & \ldots \\ 357 & 1986 & 11598 & 70125 & \ldots \\  \vdots &
\vdots & \vdots & \vdots & \ddots\end{array}\right)$$ We find that
$$L^{-1}=\left(\begin{array}{ccccc} 1 & 0 & 0 & 0 & \ldots \\0 & \frac{1}{4} & 0 & 0
& \ldots \\ 0 & 0 & \frac{1}{17} & 0 & \ldots \\ 0 & 0 & 0 &
\frac{1}{73} & \ldots
\\  \vdots & \vdots & \vdots & \vdots &
\ddots\end{array}\right)\left(\begin{array}{ccccc} 1 & 0 & 0 & 0 &
\ldots \\-15 & 4 & 0 & 0 & \ldots \\ 198 & -131 & 17 & 0 & \ldots \\
-2565 & 2875 & -854 & 73 & \ldots
\\  \vdots & \vdots & \vdots & \vdots &
\ddots\end{array}\right)$$ where the sequence $b_3(n)$ or
$1,4,17,73,314,\ldots$, \seqnum{A018902} has generating function
$\frac{1-x}{1-5x+3x^2}$ and
\begin{eqnarray*}b_3(n)&=&\sum_{k=0}^{\lfloor \frac{n}{2}
\rfloor}\binom{n-k}{k}(-3)^k5^{n-2k}-\sum_{k=0}^{\lfloor
\frac{n-1}{2} \rfloor}\binom{n-k-1}{k}(-3)^k5^{n-2k-1}\\
&=&\sum_{k=0}^n\binom{n}{k}\sum_{j=0}^k
\binom{j}{k-j}3^{2j-k}\\
&=&\sum_{k=0}^n\binom{n}{k}\sum_{j=0}^{\lfloor \frac{k}{2} \rfloor}
\binom{k-j}{j}3^{k-2j}.\end{eqnarray*} Then $3^{\binom{n}{2}}b_3(n)$
is the sequence $1, 4, 51, 1971, 228906,\ldots$. In other words, we
have
$$h_n(3)=3^{\binom{n+1}{2}}b_3(n+1).$$ We now note that $F(2n+1)$ has
generating function $\frac{1-x}{1-3x+x^2}$ with
$$F(2n+1)=\sum_{k=0}^{\lfloor \frac{n}{2} \rfloor}\binom{n-k}{k}(-1)^k3^{n-2k}-\sum_{k=0}^{\lfloor \frac{n-1}{2} \rfloor}\binom{n-k-1}{k}(-1)^k3^{n-2k-1}.$$
We can generalize this result as follows.
\begin{proposition} Let $h_n(r)$ be the Hankel transform of the sum
of the consecutive generalized Catalan numbers $c(n;r)+c(n+1;r)$.
Then $$h_n(r)=r^{\binom{n+1}{2}}(\sum_{k=0}^{\lfloor \frac{n+1}{2}
\rfloor}\binom{n-k+1}{k}(-r)^k(r+2)^{n-2k+1}-\sum_{k=0}^{\lfloor
\frac{n}{2} \rfloor}\binom{n-k}{k}(-r)^k(r+2)^{n-2k}).$$ In other
words, $h_n(r)$ is the product of $r^{\binom{n+1}{2}}$ and the
$(n+1)$-st term of the sequence with generating function
$\frac{1-x}{1-(r+2)x+rx^2}$. Equivalently, \begin{eqnarray*}
h_n(r)&=&r^{\binom{n+1}{2}}(\sum_{k=0}^{n+1}\binom{k}{n-k+1}(r+2)^{2k-n-1}(-r)^{n-k+1}-\sum_{k=0}^n\binom{k}{n-k}(r+2)^{2k-n}(-r)^{n-k})\\
&=&r^{\binom{n+1}{2}}\sum_{k=0}^{n+1}\binom{n+1}{k}\sum_{j=0}^k
\binom{j}{k-j}r^{2j-k}\\
&=&r^{\binom{n+1}{2}}\sum_{k=0}^{n+1}\binom{n+1}{k}\sum_{j=0}^{\lfloor
\frac{k}{2} \rfloor} \binom{k-j}{j}r^{k-2j}.
\end{eqnarray*}\end{proposition}
The two last expressions are a result of the fact that
$\frac{1-x}{1-(r+2)x+rx^2}$ is the binomial transform of
$\frac{1}{1-rx-x^2}$.
\section{Berlekamp Massey triangles associated to generalized Catalan numbers}
A natural question that arises when dealing with Hankel matrices is
one that is inspired by consideration of the Hankel matrix
interpretation of the Berlekamp Massey algorithm \cite{BMvar}. In
our context, this question is that of characterizing the solutions
of equations such as the following (taking $c(n,3)$ as an example)
$$\left(\begin{array}{cccc} 1 & 3 & 12
& 57  \\3 & 12 & 57 & 300  \\ 12 & 57 & 300 & 1686  \\ 57 & 300 &
1686 & 9912 \end{array}\right)
\left(\begin{array}{c}g_{1}\\g_{2}\\g_{3}\\g_{4}\end{array}\right)=\left(\begin{array}{c}300\\1686\\9912\\60213\end{array}\right)$$
or
$$\left(\begin{array}{c}g_{1}\\g_{2}\\g_{3}\\g_{4}\end{array}\right)=\left(\begin{array}{c}-81\\142\\-75\\15\end{array}\right)$$
We can define the \emph{B-M triangle} of a sequence $a_1,a_2,a_3,\ldots$ to be the lower triangular matrix whose $n$-th row is
the solution of the Berlekamp Massey equations determined by the $n$-th order Hankel matrix of the sequence.
\begin{example} The B-M triangle of the Catalan numbers. We must solve
$$(1)(x)=(1)$$
$$\left(\begin{array}{cc}1 & 1\\1 & 2\end{array}\right)\left(\begin{array}{c}x\\y\end{array}\right)
=\left(\begin{array}{c}2\\5\end{array}\right)$$
$$\left(\begin{array}{ccc}1 & 1 & 2\\1 & 2 & 5\\2 & 5 & 14\end{array}\right)\left(\begin{array}{c}x\\y\\z\end{array}\right)
=\left(\begin{array}{c}5\\14\\42\end{array}\right)$$ and so on. We obtain the triangle
$$\left(\begin{array}{ccccc} 1 & 0 & 0
& 0 & \ldots \\-1 & 3 & 0 & 0 & \ldots \\ 1 & -6 & 5 & 0 &
\ldots \\ -1 & 10 & -15 & 7 & \ldots \\  \vdots & \vdots &
\vdots & \vdots & \ddots\end{array}\right)$$
with general term $(-1)^{n-k}(\binom{n+k+1}{2k}-\binom{0}{n-k+1})$ and generating function
$$\frac{(1+x)+xy}{(1-xy)(1+2x+x^2-xy)}.$$ Regarding the entries as polynomial coefficients, we see that these polynomials are related to the
Morgan-Voyce polynomials, themselves a transformation of the Jacobi polynomials.\end{example}
In fact, the generating function of the B-M triangle for $c(n;r)$ has generating function
$$\frac{r(1+x)+xy}{(1-xy)(1+(r+1)x+rx^2-xy)}.$$
These matrices are closely related to the matrices $L$ already
studied. The standard Berlekamp Massey theory studies the
polynomials $$x^d-\sum_{i=0}^{d-1}g_ix^i.$$ For the above example,
this is $$x^4-15x^3+75x^2-142x+81.$$ Note that the companion matrix
of this polynomial is given by
$$\left(\begin{array}{cccc} 1 & 3 & 12
& 57  \\3 & 12 & 57 & 300  \\ 12 & 57 & 300 & 1686  \\ 57 & 300 &
1686 & 9912 \end{array}\right)^{-1}\left(\begin{array}{cccc} 3 & 12
& 57& 300  \\12 & 57 & 300 & 1686  \\ 57 & 300 & 1686 & 9912  \\ 300
& 1686& 9912 & 60213 \end{array}\right)=\left(\begin{array}{cccc} 0&
0 & 0& -81  \\1 & 0 & 0 & 142  \\ 0 & 1 & 0 & -75
\\ 0 & 0& 1 & 15 \end{array}\right)$$ In other words, the
characteristic polynomial of the last matrix is
$$81-142x+75x^2-15x^3+x^4.$$
Thus we get  \begin{eqnarray*}(1)^{-1}(3)&=&(3) \quad\mathrm{gives}\quad -3+x\\
\left(\begin{array}{cc}1 & 3\\3 &
12\end{array}\right)^{-1}\left(\begin{array}{c}12\\57\end{array}\right)
&=&\left(\begin{array}{c}-9\\7\end{array}\right)
\quad\mathrm{gives}\quad 9-7x+x^2\\
\left(\begin{array}{ccc}1 & 3 & 12\\3 & 12 & 57\\12 & 57 &
300\end{array}\right)^{-1}\left(\begin{array}{c}57\\300\\1686\end{array}\right)
&=&\left(\begin{array}{c}27\\-34\\11\end{array}\right)
\quad\mathrm{gives}\quad -27+34x-11x^2+x^3.\end{eqnarray*} and so
on. Forming the matrix of coefficients, we obtain
$$\left(\begin{array}{cccccc} 1 & 0 & 0 & 0 & 0 &\ldots \\-3 & 1 & 0
& 0& 0 &\ldots \\ 9 & -7 & 1 & 0 & 0&\ldots \\ -27 & 34 & -11 & 1
&0& \ldots
\\ 81 & -142 & 75 & -15 & 1 \\ \vdots & \vdots & \vdots & \vdots &  \vdots & \ddots\end{array}\right)$$
which is the Riordan matrix
$$(\frac{1}{1+3x},\frac{x}{1+4x+3x^2}).$$
The inverse of this matrix is the $L$ matrix in the $LDL^T$
decomposition of the Hankel matrix for $c(n;3)$. The corresponding
B-M matrix as defined above is given by the negative of the
sub-diagonal matrix.

\bigskip
\hrule
\bigskip
\noindent 2000 {\it Mathematics Subject Classification}: Primary
11B83; Secondary 05E35,11C20,11Y55,15A23.

\noindent \emph{Keywords:} Hankel transform, Fibonacci numbers,
Catalan numbers, Riordan array, Orthogonal polynomials.

\bigskip
\hrule
\bigskip
Concerns sequences  \seqnum{A000108}, \seqnum{A000129},
\seqnum{A001333}, \seqnum{A007318}, \seqnum{A008288},
\seqnum{A018902}, \seqnum{A079496}

\end{document}